\documentclass[12pt,a4paper]{amsart}
\usepackage{amssymb,amsmath}

\headsep=1cm \numberwithin{equation}{section}
\newtheorem{theorem}{Theorem}[section]
\newtheorem{lemma}[theorem]{Lemma}
\newtheorem{proposition}[theorem]{Proposition}
\newtheorem{corollary}[theorem]{Corollary}

\theoremstyle{definition}
\newtheorem{definition}[theorem]{Definition}
\theoremstyle{remark}

\newtheorem{example}[theorem]{Example}

\newcommand{\Ass}{\operatorname{Ass}}

\newcommand{\Spec}{\operatorname{Spec}}

\newcommand{\rad}{\operatorname{rad}}

\newcommand{\Supp}{\operatorname{Supp}}

\newcommand{\Ann}{\operatorname{Ann}}

\newcommand{\Max}{\operatorname{Max}}

\newcommand{\fp}{\frak{p}}
\newcommand{\fq}{\frak{q}}
\newcommand{\fa}{\frak{a}}
\newcommand{\fb}{\frak{b}}

\begin{document}
\author[Divaani-Aazar and Esmkhani ]{Kamran Divaani-Aazar and Mohammad Ali
Esmkhani}
\title[Associated prime submodules of finitely generated
modules]{Associated prime submodules of finitely generated
modules}

\address{K. Divaani-Aazar, Department of Mathematics, Az-Zahra University,
Vanak, Post Code 19834, Tehran, Iran and Institute for Studies in
Theoretical Physics and Mathematics, P.O. Box 19395-5746, Tehran,
Iran.} \email{kdivaani@ipm.ir}

\address{M.A. Esmkhani, Department of Mathematics, Shahid Beheshti University,
Tehran, Iran and Institute for Studies in Theoretical Physics and
Mathematics, P.O. Box 19395-5746, Tehran, Iran.}



\begin{abstract}
Let $R$ be  a commutative ring with identity. For a finitely
generated $R$-module $M$, the notion of associated prime
submodules of $M$ is defined. It is shown that this notion
inherits most of essential properties of the usual notion of
associated prime ideals. In particular, it is proved that for a
Noetherian multiplication module $M$, the set of associated prime
submodules of $M$ coincides with the set of $M$-radicals of
primary submodules of $M$ which appear in a minimal primary
decomposition of the zero submodule of $M$. Also, Anderson's
theorem [{\bf 2}] is extended to minimal prime submodules in a
certain type of modules.
\end{abstract}

\maketitle

\section{Introduction}

Recently, extensive research has been done  on prime submodules.
Let $R$ be a commutative ring with identity and $M$ an $R$-module.
A proper submodule $N$ of $M$ is said to be prime (or $\fp$-prime)
if $re\in N$ for $r\in R$ and $e\in M$, implies that either $e\in
N$ or $r\in \fp=N:M$. A general theme in the studying of prime
submodules is to extend results concerning prime ideals to prime
submodules. For example, Cohen's Theorem, Prime Avoidance Theorem
and Krull's Principal Ideal Theorem are generalized to prime
submodules in [{\bf 13}, Theorem 5 ], [{\bf 8}, Theorem 2.3] and
[{\bf 6}, Theorem 11] respectively.

The use of the notion of  associated prime ideals has found
substantial applications in commutative algebra. In fact, the set
of associated prime ideals of a module contains a lot of
information about the module itself. It is natural to expect that,
if there is an appropriate definition of associated prime
submodules, then many results concerning associated prime ideals
can be generalized to associated prime submodules.

The main goal of this paper is to define the concept of associated
prime submodules of finitely generated modules and to investigate
their properties. Surely, there is much more work to be done. Let
$M$ be a finitely generated $R$-module and $\fp$ a prime ideal of
$R$. Following [{\bf 11}], we set
$$M(\fp)=\{x\in M: sx\in \fp M,\  \text{for\ some}\ s\in R\smallsetminus\fp\}.$$
It is  easy to see that, if $\fp$ contains the annihilator of $M$,
then $M(\fp)$ is a $\fp$-prime submodule of $M$. We define the set
of associated prime submodules of $M$ as $$\Ass_PM=\{M(\fp): \fp
\in \Ass_RM\},$$ where $\Ass_RM$ denotes the set of associated
prime ideals of $M$.

In section 2, some results related to associated prime ideals are
extended to associated prime submodules. In particular, it is
shown that, if $0=\bigcap^n_{i=1}Q_i$ is a minimal primary
decomposition of the zero submodule of the  Noetherian
multiplication module M, then  $\Ass_PM=\{\rad(Q_i): i=1,2,\dots
,n\}.$ Recall that an $R$-module $M$ is called multiplication, if
every submodule $N$ of $M$ is of the form $\fa M$, for some ideal
$\fa$ of $R$. Also, for a submodule $N$ of the $R$-module $M$,
$\rad(N)$, the $M$-radical of $N$, is defined as the intersection
of the prime submodules of $M$ containing $N$ (see [{\bf 10}]).

In section 3, we study the set of minimal associated prime
submodules of finitely generated modules. Firstly, we introduce
the class of quasi multiplication modules. This class of modules
contains multiplication, finitely generated weak multiplication
and flat modules. Then, we extend Anderson's theorem to minimal
prime submodules of finitely generated quasi multiplication
modules to the effect that if all minimal prime submodules of a
finitely generated quasi multiplication $R$-module $M$ are
finitely generated, then their number is finite.

Throughout this paper, $R$ is a commutative ring  with identity
and all modules are assumed to be unitary.

\section{Associated prime submodules}

Recall that the sets of associated and of supported prime ideals
of a given $R$-module $M$ are defined respectively as:
$$\Ass_RM=\{\fp \in \Spec R: \fp =(0:x), \text{for\ some\ nonzero\
element}\ x \ \text {of}\ R \},$$and
$$\Supp_RM=\{\fp \in \Spec R:
\fp \supseteq (0:x), \text{for\ some\ nonzero\ element}\ x \ \text
{of}\ R \}.$$ It is well known that, if $M$ is finitely generated,
then $\Supp_RM$ is the set of all prime ideals of $R$ which
contain $\Ann_RM$.

For completeness, we collect some important known properties of
these notions in the following lemma.

\begin{lemma} i) If $R$ is Noetherian, then $M$ is zero if and
only if $\Ass_RM$ is empty.\\
ii) If $R$ is Noetherian, then the set of minimal elements of
 $\Ass_RM$ and that of $\Supp_RM$ are equal.\\
iii) If $M$ is Noetherian, then $\Ass_RM$ is finite.\\
iv) If $M$ is Artinian, then $\Ass_RM=\Supp_RM$,  and this set
consists of finitely many maximal ideals. \\
v) If $R$ is Noetherian and $S$ a multiplicatively closed subset
of $R$, then $$\Ass_{S^{-1}R}S^{-1}M=\{S^{-1}\fp: \fp \in \Ass_RM
\text{and}\ \fp \cap S=\phi\}.$$ vi) If $R$ is Noetherian and
$0=\bigcap^n_{i=1}Q_i$ is a minimal primary decomposition of the
zero submodule of $M$, then $\Ass_RM=\{\rad (Q_i:M): i=1,2,\dots
,n\}$.
\end{lemma}

{\bf Proof.} (i) and (ii) hold respectively, by  [{\bf 9}, 7.B
Corollary 1] and [{\bf 9}, Theorem 9]. Also, (v) and (vi) follow
respectively, by [{\bf 9}, 7.C Lemma] and [{\bf 9}, 8.E Lemma]. It
is easy to deduce (iv) from [{\bf 12}, Exercise 8.49]. Finally,
(iii) is clear by (vi). $\Box$

In the sequel, we generalize the above mentioned properties of
associated prime ideals to associated prime submodules. Let $M$ be
an $R$-module and $\fp$ a prime ideal of $R$. Following [{\bf
11}], we denote the set $\{x\in M: sx\in \fp M,\  for\ some\ s\in
R\smallsetminus\fp\}$, by $M(\fp)$. We summarize some important
properties of this notion in the following lemma and in the
sequel, we may use them without further comment.

\begin{lemma} Let $M$ be an $R$-module and $\fp$ a prime ideal of $R$.\\
i) $M=M(\fp)$ or $M(\fp)$ is a $\fp$-prime submodule of $M$.\\
ii) If $M$ is either finitely generated or multiplication, then
$M(\fp)$ is a $\fp$-prime submodule of $M$ if and only if
$\fp\in \Supp_RM$.\\
iii) Every $\fp$-prime submodule of $M$ contains $M(\fp)$.
\end{lemma}

{\bf Proof.} (i) holds by [{\bf 11}, Proposition 1.7], while (iii)
is true by [{\bf 11}, Lemma 1.6]. \\
Now, we show (ii). Suppose $M$ is either finitely generated or
multiplication.  It is easy to check that  $M=M(\fp)$ if and only
if $M_{\fp}=(\fp R_{\fp})M_{\fp}$. By [{\bf 3}, Proposition 1],
the assertion of Nakayama's Lemma holds also for multiplication
modules. Thus, in both cases, it follows that $M_{\fp}=(\fp
R_{\fp})M_{\fp}$ if and only if $M_{\fp}=0$. Therefore, it turns
out by (i), that $M(\fp)$ is a $\fp$-prime submodule of $M$ if and
only if $\fp\in \Supp_RM$. $\Box$

Now, we are ready to present the definitions of associated and of
supported prime submodules of a finitely generated module.

\begin{definition} (i) Let $M$ be an $R$-module. We say $M$ is
{\it weakly finitely generated}, if for any $\fp\in \Supp_RM$, the
submodule $M(\fp)$ of $M$ is proper.\\
(ii) Let $M$ be a weakly finitely generated $R$-module. We define
the sets of associated and of supported prime submodules of $M$,
respectively as:$$\Ass_PM=\{M(\fp): \fp\in \Ass_RM\},$$ and
$$\Supp_PM=\{M(\fp): \fp\in \Supp_RM\}.$$
\end{definition}

\begin{example} By Lemma 2.2(ii), it becomes clear that if the
$R$-module $M$ is either finitely generated or multiplication,
then $M$ is weakly finitely generated. It is worthy to mention
that, by [{\bf 5}, Corollary 3.9] over a Noetherian ring, any
multiplication module is finitely generated.
\end{example}

If $M$ is a Noetherian $R$-module, then the ring $T=R/\Ann_RM$ is
a Noetherian ring. One can check easily that $M$ is a $T$-module
and that $\Ass_TM=\{\fp/\Ann_RM: \fp\in \Ass_RM \}$. Hence the
following is immediate, by Lemma 2.1.

\begin{lemma} Suppose that $M$ is a Noetherian $R$-module.\\
i) $\Ass_PM$ is finite.\\
ii) $M=0$ if and only if $\Ass_PM=\phi$.
\end{lemma}

The set of all prime submodules of the $R$-module $M$ is denoted
by $\Spec M$. Also, the set of maximal submodules of $M$ is
denoted by $\Max M$.

\begin{lemma} Assume that $M$ is an $R$-module of finite length.\\
i) $\Ass_PM=\Supp_PM$ and this set is finite.\\
ii) Moreover, if $M$ is a multiplication module, then
$\Ass_PM=\Max M=\Spec M$.
\end{lemma}

{\bf Proof.} i) follows by Lemma 2.1(iv).\\
ii) Let $M(\fp)$ be an associated prime submodules of $M$. Then,
by Lemma 2.1(iv), $\fp$ is a maximal ideal of $R$. Let $N$ be a
proper submodule of $M$ containing $M(\fp)$. Then
$$\fp=(M(\fp):M)\subseteq (N:M)\subsetneqq R.$$ Hence $\fp=(N:M)$.
But $M$ is a multiplication module and so
$$N=(N:M)M=(M(\fp):M)M=M(\fp).$$ Thus $M(\fp)$ is a maximal
submodule of $M$.

Because any maximal submodule of $M$ is a prime submodule, to
complete the proof, it suffices to show that any prime submodule
of $M$ is an associated prime submodule. Assume that $N$ is a
prime submodule of $M$. It follows that $\fp=(N:M)$ is an element
of $\Supp_RM$ and so $\fp$ belongs to $\Ass_RM$, by Lemma 2.1(iv).
Now, we have $N= \fp M = M(\fp)$, and so $N\in \Ass_PM$, as
required. $\Box$

\begin{corollary} Let $M$ be a  multiplication module. The following
are equivalent:\\
i) $M$ is Artinian.\\
ii) $M$ is Noetherian and $\Ass_PM\subseteq \Max M$.
\end{corollary}

{\bf Proof.} By [{\bf 5}, Corollary 2.9], any Artinian
multiplication module is cyclic.  Hence (i) implies (ii), by Lemma
2.6, and the fact that every finitely generated Artinian module is
Noetherian.

Now, assume (ii) holds. Let $\fp\in \Ass_RM$. Then $M(\fp)\in \Max
M$, by the assumption. This yields that $\fp\in \Max R$. Using
Lemma 2.1(ii), we can deduce that every prime ideal of the
Noetherian ring $T=R/\Ann_RM$ is maximal, and so T is an Artinian
ring. Now, because $M$ is finitely generated, we can deduce that
$M$ is Artinian as an $R$-module. $\Box$

\begin{lemma} Let $S$ be a multiplicatively closed subset of the
Noetherian ring $R$ and let $M$ be a finitely generated
$R$-module. Then the set of associated prime submodules of the
$S^{-1}R$-module $S^{-1}M$ is equal to $\{S^{-1}P: P\in \Ass_PM \
\text{and}\ (P:M)\cap S=\phi\}.$
\end{lemma}

{\bf Proof.} By Lemma 2.1(v), $\Ass_{S^{-1}R}S^{-1}M=\{S^{-1}\fp:
\fp\in \Ass_RM  \  \text{and}\ \fp\cap S=\phi\}$. Let $P\in
\Ass_PM$ be such that $(P:M)\cap S=\phi$. Set $\fp=P:M$. Then
$\fp$ is an associated prime ideal of $M$ and $P=M(\fp)$. It is
straightforward to see that $S^{-1}P=S^{-1}M(S^{-1}\fp)$, and so
$S^{-1}P$ is an associated prime submodule of the $S^{-1}R$-module
$S^{-1}M$.

Conversely, assume that $Q$ is an associated prime submodule of
$S^{-1}M$ as an $S^{-1}R$-module. Then there exists $\fp\in
\Ass_RM$, with $\fp \cap S=\phi$, such that
$Q=S^{-1}M(S^{-1}\fp)$. Thus $Q=S^{-1}(M(\fp))$ and $M(\fp)\in
\Ass_PM$. $\Box$

Now, we are ready to prove the main result of this section.

\begin{theorem} Let $M$ be a Noetherian multiplication $R$-module.
If $0=\bigcap^n_{i=1}Q_i$ is a minimal primary decomposition of
the zero submodule of $M$, then $\Ass_PM=\{\rad(Q_i): i=1,2,\dots
,n\}$.
\end{theorem}

{\bf Proof.} Set $\fp_i=\rad (Q_i:M)$, for $i=1,2,\dots ,n$. Then
by Lemma 2.1(vi), $\Ass_RM=\{\fp_i: i=1,2,\dots ,n\}$. Fix $1\leq
i\leq n$. It turns out by [{\bf 10}, Theorem 4], that
$\rad(Q_i)=(\rad (Q_i:M))M$. Now, we have
$$M(\fp_i)=(M(\fp_i):M)M=\fp_iM.$$
Hence $\rad(Q_i)\in \Ass_PM$.

Conversely, assume that $P\in \Ass_PM$. Then there is $1\leq i\leq
n$, such that $P=M(\fp_i)$. But the assumption on $M$ implies that
$M(\fp_i)=\fp_iM$. Thus
$$P=(\rad (Q_i:M))M=\rad(Q_i),$$ as required. $\Box$

\section{Minimal associated prime submodules}

An $R$-module $M$ is called a weak multiplication module if every
prime submodule $P$ of $M$ is of the form $\fp M$, for some prime
ideal $\fp$ of $R$ [{\bf 1}]. Next, we present the following
definition.

\begin{definition} An $R$-module $M$ is called a {\it quasi
multiplication module} if $M(\fp)=\fp M$, for all $\fp\in
\Supp_RM$.
\end{definition}

\begin{example} i) Let $M$ be a weakly finitely generated $R$-module
which is weak multiplication. Then $M$ is quasi multiplication.\\
ii) Every flat $R$-module is a quasi multiplication module. To see
this, let $\fp\in \Supp_RM$. If $\fp M=M$, then $M(\fp)=\fp M$,
because $\fp M\subseteq M(\fp)$, as one can see clearly. Now
assume that $\fp M$ is a proper submodule of $M$. Then $\fp M$ is
a $\fp$-prime
submodule of $M$, by [{\bf 7}, Theorem 3]. Hence $M(\fp)=\fp M$,
by Lemma 2.2(iii)\\
iii) It turns out, by Lemma 2.2(ii), that every multiplication
$R$-module is quasi multiplication.
\end{example}

The following is clear, by definition of quasi multiplication
modules.

\begin{lemma} Let $M$ be a weakly finitely generated module. If $M$
 is quasi multiplication, then $\Ass_PM=\{\fp M: \fp\in \Ass_RM \}$
 and $\Supp_PM=\{\fp M:\fp\in \Supp_RM \}$.
\end{lemma}

\begin{proposition} Let $M$ be a finitely generated $R$-module.
Assume that $M$ is quasi multiplication. Then the following hold.\\
i) The set of minimal prime submodule of $M$ is equal to
$$\{\fp M:
\fp \  \text {is\ a\ minimal\ element\ of}\  \Supp_RM \}.$$ ii) If
$R$ is Noetherian, then the set of minimal elements of the $\Spec
M$, $\Supp_PM$ and $\Ass_PM$ are coincide. Consequently, the set
of minimal elements of $\Spec M$, $\Supp_PM$ and $\Ass_PM$ all
coincide, whenever $M$ is a Noetherian module over an arbitrary
ring.
\end{proposition}

{\bf Proof.} i) Let $\fp$ be a minimal element of $\Supp_RM$.
Since $M$ is  a finitely generated quasi multiplication module, it
follows that $M(\fp)=\fp M$, and so $\fp M$ is a prime submodule
of $M$. Assume that $Q$ is a prime submodule of $M$ such that
$Q\subseteq \fp M$. Then
$$\Ann_RM\subseteq Q:M\subseteq \fp M:M=\fp.$$ Hence $Q:M=\fp$, by the assumption on
$\fp$. It turns out by Lemma 2.2(iii), that $Q=M(\fp)=\fp M$.

Now, assume that $P$ is a minimal prime submodule of $M$. Let
$\fp=P:M$. Then $P=M(\fp)$, by Lemma 2.2(iii). Let $\fq$ be an
element of $\Supp_RM$ such that $\fq\subseteq \fp$. Then, it
follows that
$$M(\fq)=\fq M\subseteq \fp M=M(\fp)=P,$$ and so $M(\fq)=P$. Thus
$\fp=P:M=M(\fq):M=\fq$, and so $\fp$ is minimal in $\Supp_RM$.\\

ii) Since $\Ass_PM\subseteq \Supp_PM\subseteq \Spec M$, it is
enough to show that every minimal associated prime submodule of
$M$ is minimal in $\Spec M$. Note that in view of part (i) and
Lemma 2.1(ii), it follows that any prime submodule of $M$ contains
an element of $\Ass_PM$. By Lemma 2.1(ii), the set of minimal
elements of $\Ass_RM$ and that of $\Supp_RM$ are equal. Let
$M(\fp)$ be a minimal element of $\Ass_PM$ and $Q$ a minimal prime
submodule of $M$, which is contained in $M(\fp)$. Then by (i),
$Q=\fq M$, for some minimal element $\fq$ of $\Supp_RM$. Hence
$\fq\in \Ass_RM$, and so $Q\in \Ass_PM$. Therefore, $Q=M(\fp)$, as
required.

Set $T=R/\Ann_RM$. Since $M$ possesses the structure of a
$T$-module in a natural way, the last assertion of part (ii)
follows immediately. $\Box$

\begin{lemma} Let $M$ be a weakly finitely generated $R$-module and let
$\{\fp_1,\dots ,\fp_n\}$ be a subset of minimal elements of
$\Supp_RM$. If $\fp_1\dots \fp_nM=0$, then $\fp_1,\dots ,\fp_n$
are the only minimal elements of $\Supp_RM$.
\end{lemma}

{\bf Proof.} Let $\fp$ be a minimal element of $\Supp_RM$. Then
$M(\fp)$ is a $\fp$-prime submodule of $M$. Thus
$$\fp_1\dots \fp_n\subseteq M(\fp):M=\fp,$$ and so $\fp=\fp_i$,
for some $1\leq i \leq n$. $\Box$

\begin{proposition} Let $M$ be a finitely generated $R$-module.
Assume that $\fp M$ is finitely generated for all minimal elements
$\fp$ of $\Supp_RM$. Then the number of minimal elements of
$\Supp_RM$ is finite.
\end{proposition}

{\bf Proof.} In view of Lemma 3.5, it suffices to show that there
are minimal elements $\fp_1,\dots ,\fp_n$ of $\Supp_RM$, such that
$\fp_1\dots \fp_nM=0$. Suppose that the contrary is true and we
search for a contradiction. Let $\mathcal{A}$ denote the set of
all ideals $\fa$ of $R$, such that $\Ann_RM\subseteq \fa$ and $\fa
M$ does not contain any submodule of the form $\fp_1\dots \fp_nM$,
where $\fp_i$'s are minimal elements of $\Supp_RM$. Then
$\mathcal{A}$ is not empty and by Zorn's lemma, we deduce that it
has a maximal element $\fp$, say. Note that, if $\fa$ and $\fb$
are two ideals of $R$ such that $\fa M$ and $\fb M$ are finitely
generated, then the submodule $\fa \fb M$ is also finitely
generated.

We show that $\fp$ is a prime ideal. To this end, let $x$ and $y$
be two elements of $R$ such that $x\notin \fp ,y\notin \fp$ and
$xy\in \fp$. Then $\fp +Rx$ and $\fp +Ry$ are not in
$\mathcal{A}$. Hence there are minimal elements $\fp_1,\dots
,\fp_n$ and $\fq_1,\dots ,\fq_m$ of $\Supp_RM$ such that
$\fp_1\dots \fp_nM\subseteq (\fp +Rx)M$ and $\fq_1\dots
\fq_mM\subseteq (\fp+Ry)M$. Then
$$(\fp_1\dots \fp_n \fq_1\dots \fq_m)M\subseteq \fp M,$$
which is a contradiction. Therefore $\fp$ is a prime ideal of $R$.
Since $\fp\in \Supp_RM$, it follows that $\fp$ contains a minimal
element $\fq$ of $\Supp_RM$, and so $\fq M\subseteq \fp M$. This
contradicts the assumption that $\fp$ belongs to $\mathcal{A}$.
$\Box$

A result due to Anderson [{\bf 2}], asserts that, if every minimal
prime ideal of $R$ is finitely generated, then the number of
minimal prime ideals of $R$ is finite. Then Anderson's result is
generalized to prime submodules in multiplication modules by
Behboodi and Koohy in [{\bf 4}, Theorem 2]. They proved that if
every minimal prime submodule of the multiplication $R$-module $M$
is finitely generated, then $M$ has only finitely many minimal
prime submodules. In view of Proposition 3.4(i) and Proposition
3.6, we can extend Anderson's theorem to finitely generated quasi
multiplication modules. Namely:

\begin{theorem} Let $M$ be a finitely generated quasi
multiplication $R$-module. Assume that every minimal prime
submodule of $M$ is finitely generated. Then $M$ has only finitely
many minimal prime submodules.
\end{theorem}



\begin{thebibliography}{99}
\bibitem{} S. Abu-Saymeh, {\it On dimensions of finietly generated
modules}, Comm. Algebra, {\bf 23}(3) (1995), 1131-1144.
\bibitem{} D. D. Anderson, {\it A note on minimal prime ideals},
Proc. AMS., {\bf 122} (1994), 13-14.
\bibitem{} A. Barnard, {\it Multiplication modules}, J. Algebra,
{\bf 71} (1981), 174-178.
\bibitem{} M. Behboodi and H. Koohy, {\it On minimal prime submodules},
Far East J. Math. Sci. (FJMS), {\bf 6}(1) (2002), 83-88.
\bibitem{} Z. A. EL-Bast and P. F. Smith, {\it Multiplication
modules}, Comm. Algebra, {\bf 16}(4) (1988), 755-779.
\bibitem{} S. M. George, R. L. McCasland and P. F. Smith, {\it A principal ideal theorem
analogue for modules over commutative rings}, Comm. Algebra, {\bf
22}(6) (1994), 2083-2099.
\bibitem{} C. P. Lu, {\it Prime submodules of modules}, Comment.
Math. Univ. St. Paul., {\bf 33}(1) (1984), 61-69.
\bibitem{} C. P. Lu, {\it Union of prime submodules}, Hoston Journal of Mathematics,
{\bf 23}(2) (1997), 203-213.
\bibitem{} H. Matsumura, {\it Commutative algebra}, Benjamin/Cummings Publishing Company.,
Inc., Reading, MA., 1980.
\bibitem{} R. L. McCasland and M. E. Moore, {\it On radical of submodules of finitely
generated modules}, Cand. Math. Bull., {\bf 29}(2) (1986), 37-39.
\bibitem{} R. L. McCasland and P. F. Smith, {\it Prime submodules of Noetherian modules},
Rocky Mountain J. Math., {\bf 23}(3) (1993), 1041--1062.
\bibitem{} R. Y. Sharp, {\it Steps in commutative algebra}, Cambridge University
Press, Cambridge 1990.
\bibitem{} P. F. Smith, {\it Concerning a theorem of I. S. Cohen}, An. Stiin. Univ.
" Ovidius" Constanta, Ser. Math., {\bf 2} (1994), 160-167.


\end{thebibliography}
\end{document}